\documentclass[11pt]{article}
\usepackage{amsmath,amsfonts,amssymb,theorem}

 \newcommand{\doubleheaddownarrow}{\big\downarrow\kern-3.325mm\downarrow}
 \newcommand{\wY}{\widetilde Y}
 \newcommand{\wH}{\widetilde H}
 \newcommand{\wD}{\widetilde D}
 \newcommand{\broken}{\dasharrow}
 \newcommand{\1}{^{-1}}
 \newcommand{\iso}{\cong}
 \newcommand{\tensor}{\otimes}
 \newcommand{\lodot}{_{\scriptscriptstyle{\bullet}}}
 \newcommand{\Fbar}{\overline F}
 \newcommand{\sbar}{\overline s}
 \newcommand{\tbar}{\overline t}
 \newcommand{\xbar}{\overline x}
 \newcommand{\fm}{\mathfrak m}
 \newcommand{\fn}{\mathfrak n}
 \newcommand{\fp}{\mathfrak p}
 \newcommand{\into}{\hookrightarrow}
 \newcommand{\onto}{\twoheadrightarrow}
 \newcommand{\rest}[1]{{}_{|#1}}
 \newcommand{\aff}{\mathbb A}
 
 \newcommand{\PP}{\mathbb P}
 \newcommand{\Q}{\mathbb Q}
 \newcommand{\Z}{\mathbb Z}
 \newcommand{\Oh}{\mathcal O}

 \newcommand{\al}{\alpha}
 
 \newcommand{\fie}{\varphi}
 
 \newcommand{\si}{\sigma}
 \newcommand{\om}{\omega}

 \renewcommand{\div}{\operatorname{div}}
 \newcommand{\coker}{\operatorname{coker}}
 \newcommand{\res}{\operatorname{res}}

 \newcommand{\Ass}{\operatorname{Ass}}
 \newcommand{\Ann}{\operatorname{Ann}}
 \newcommand{\depth}{\operatorname{depth}}
 \newcommand{\codim}{\operatorname{codim}}
 \newcommand{\socle}{\operatorname{socle}}
 
 \newcommand{\Ext}{\operatorname{Ext}}
 \newcommand{\sExt}{\operatorname{\mathcal E\mathit{xt}}}
 \newcommand{\Proj}{\operatorname{Proj}}
 \newcommand{\sHom}{\operatorname{\mathcal H\mathit{om}}}
 \newcommand{\Hom}{\operatorname{Hom}}
 
 \newcommand{\Sing}{\operatorname{Sing}}
 \newcommand{\Spec}{\operatorname{Spec}}

 \newtheorem{theorem}{Theorem}[section]
 \newtheorem{lemma}[theorem]{Lemma}
 \newtheorem{prop}[theorem]{Proposition}
 
 {
 \theorembodyfont{\rmfamily}
 \newtheorem{defn}[theorem]{Definition}
 
 \newtheorem{exc}[theorem]{Exercise}
 \newtheorem{ex}[subsection]{Exercise}
 \newtheorem{rem}[theorem]{Remark}
 
 \newtheorem{rmk}[subsection]{Remark}

 }
 \newenvironment{pf}{\paragraph{Proof}}{\par\medskip}
 
 \newcommand{\qed}{\ifhmode\unskip\nobreak\fi\quad\ensuremath\square}
\newcommand{\QED}{\ifhmode\unskip\nobreak\fi\quad\ensuremath{\mathrm{QED}}}

\numberwithin{equation}{section}

\title{Kustin--Miller unprojection without complexes}

\author{Stavros Papadakis\thanks{The first author thanks
the Greek State Scholarships Foundation for support.}
\and Miles Reid\thanks{We both thank Kyoto University, RIMS for
generous support and hospitality.}}
\date{Nov 2001}

\begin{document}
\maketitle
 \begin{abstract}
 Gorenstein projection plays a key role in birational geometry; the
typical example is the linear projection of a del Pezzo surface of
degree $d$ to one of degree~$d-1$, but variations on the same idea
provide many of the classical and modern birational links between Fano
3-folds. The inverse operation is the Kustin--Miller unprojection
theorem, which constructs ``more complicated'' Gorenstein rings
starting from ``less complicated'' ones (increasing the codimension by
1). We give a clean statement and proof of their theorem, using the
adjunction formula for the dualising sheaf in place of their complexes
and Buchsbaum--Eisenbud exactness criterion. Our methods are scheme
theoretic and work without any mention of the ambient space. They are
thus not restricted to the local situation, and are well adapted to
generalisations.

Section~\ref{sec:appl} contains examples, and discusses briefly the
applications to graded rings and birational geometry that motivate this
study; see also Papadakis \cite{P1} and Reid \cite{Ki}--\cite{T4}.
 \end{abstract}

\section{The theorem}

Let $X=\Spec\Oh_X$ be a Gorenstein local scheme and $I\subset\Oh_X$ an
ideal sheaf defining a subscheme $D=V(I)\subset X$ that is also Gorenstein
and has codimension~1 in $X$. We assume that all schemes are Noetherian. We
do not assume anything else about the singularities of $X$ and $D$,
although an important case in applications is when $X$ is normal and $D$ a
Weil divisor.

The adjunction formula (compare Reid \cite{R}, Appendix to Section~2) gives
 \[
 \om_D=\sExt^1(\Oh_D,\om_X).
 \]
To calculate the $\sExt$, we apply the derived functor of $\sHom$ to the
exact sequence $0\to I\to\Oh_X\to\Oh_D\to0$ into $\om_X$, obtaining the
usual adjunction exact sequence
 \[
 0\to\om_X\to\sHom(I,\om_X)\xrightarrow{\,\res_D\,}\om_D\to0,
 \]
where $\res_D$ is the residue map. For example, in the case that $X$ is
normal and $D$ a divisor, the second map is the standard Poincar\'e
residue map $\Oh_X(K_X+D)\to\Oh_D(K_D)$.

\begin{lemma}\label{lem!inj} The $\Oh_X$ module $\sHom(I,\om_X)$ is
generated by two elements $i$ and $s$, where $i$ is a basis of\/ $\om_X$
and $s\in\sHom(I,\om_X)$ satisfies
 \begin{enumerate}
 \renewcommand{\labelenumi}{(\roman{enumi})}
 \item $s\colon I\into\om_X$ is injective;
 \item $\sbar=\res_D(s)$ is a basis of $\om_D$.
 \end{enumerate}
\end{lemma}

\begin{pf} Choose bases $i\in\om_X$, $\sbar\in\om_D$ and any lift
$s\mapsto\sbar$. Then everything holds except (i). We achieve (i) by a
simple exercise in primary decomposition: write $X_i\subset X$ for the
reduced irreducible components and $P_i$ for the corresponding minimal
prime ideals, so that $X_i=V(P_i)$. Then $\Ass\Oh_X=\{P_i\}$ because $X$
is Cohen--Macaulay.

We have $\ker s\subset I$, so its Ass also consists of irreducible
components, namely those on which $s$ vanishes. Choose $f\in\Oh_X$ such
that
 \[
 \begin{cases}
 f\notin P_i & \hbox{for } P_i\in\Ass(\ker s), \\
 f\in P_i & \hbox{for } P_i\notin\Ass(\ker s);
 \end{cases}
 \]
in other words, $f$ is nonzero (thus generically a unit) on each component
where $s$ vanishes, and $f$ vanishes along every component at which $s$ is
a unit. Now replacing $s$ by $s+fi$ gives (i). \QED \end{pf}

We view $s$ as defining an isomorphism $I\to J$, where $J\subset\om_X=
\Oh_X$ is another ideal. Choose a set of generators $f_1,\dots,f_k$ of
$I$ and write $s(f_i)=g_i$ for the corresponding generators of $J$. We
view $s=g_i/f_i$ as a rational function having $I$ as ideal of
denominators and $J$ as ideal of numerators (compare Remark~\ref{rmk!expl}).
Unprojection is simply the graph of $s$.

 \begin{defn}\label{defn!un}
 Let $S$ be an indeterminate. The {\em unprojection ring} of $D$ in $X$ is
the ring $\Oh_X[s]=\Oh_X[S]/(Sf_i-g_i)$; the {\em unprojection} of $D$ in
$X$ is its Spec, that is,
 \[
 Y=\Spec\Oh_X[s].
 \]
Clearly, $Y$ is simply the subscheme of $\Spec\Oh_X[S]=\aff^1_X$ defined
by the ideal $(Sf_i-g_i)$. Usually $Y$ is no longer local (see
Example~\ref{ssec!nodal}).
 \end{defn}

\begin{rem} \label{rmk!expl}
 \begin{enumerate}
 \renewcommand{\labelenumi}{(\arabic{enumi})}
 \item Clearly $J=\Oh_X$ if and only if $I$ is principal; if $I=(f)$
then $\Oh_X[s]=\Oh_X[1/f]$. We exclude this case in what follows.

 \item We only choose generators for ease of notation here. The ideal
defining $Y$ could be written $\bigl(Sf-s(f)\bigm|f\in I\bigr)$.

 The construction is independent of $s$: the only choice in
Lemma~\ref{lem!inj} is $s\mapsto us+hi$ with $u,h\in\Oh_X$ and $u$ a unit,
which just gives the affine linear coordinate change $S\mapsto uS+h$ in
$\aff^1_X$.

 \item \label{rmk!ratio} The total ring of fractions $K(X)$ is defined
as $S\1\Oh_X$ where $S$ is the set of non-zerodivisors, that is, the
complement of the associated primes $P_i\in\Ass\Oh_X$. Then $s\colon I\to
J$ is multiplication by an invertible rational function in $K(X)$. For $I$
contains a regular element $w$ (in fact $ \depth I=\codim D=1$, by
Matsumura \cite{M}, Theorem~17.4), and $s(w)/w \in K(X)$ is independent of
the choice of $w$, because
 \[
 0=s(w_1w_2-w_2w_1)=w_1s(w_2)-w_2s(w_1) \quad \hbox{ for } 
 w_1,w_2 \in I.
 \]

 \item We defined $\Oh_X[s]$ by generators and relations in
Definition~\ref{defn!un}. If $X$ is normal, it equals the subring of
$K(X)$ generated by $\Oh_X$ and $s$. Sketch proof: The point is to prove
``only linear relations'', that is, any relation of the form
$as^2+bs+c=0$ (etc.)\ is in the ideal generated by the linear relations
$sf_i-g_i$; this is clear, because if $(as+b)s=-c\in\Oh_X$ then $as+b=-c/s$
cannot have any divisor of poles.
 \end{enumerate}
\end{rem}

\begin{lemma} Write $N=V(J)\subset X$ for the subscheme with $\Oh_N=\Oh_X/J$.
\begin{enumerate}
 \item[(a)] No component of $X$ is contained in $N$.
 \item[(b)] Every associated prime of\/ $\Oh_N$ has codimension~$1$.
\end{enumerate}
\end{lemma}

If $X$ is normal then $D$ and $N$ are both divisors, with $\div s=N-D$. More
generally, set $n=\dim X$; then (a) says that $\dim N\le n-1$, and (b) says
that $\dim N=n-1$ (and has no embedded primes).

\begin{pf} As we have just said, $I$ contains a regular element
$w\in\Oh_X$. Then $v=s(w)\in J$ is again regular (obvious), and (a)
follows.

(b) follows from (a) and Step~2 of the proof of Theorem~\ref{thm!KM}: for
(\ref{eq!dpth}) below gives $n-1\le \depth\Oh_X/J<n=\dim X$.

For a direct proof of (b), note first that $vI=wJ$. We prove that every
element of $\Ass(\Oh_X/vI)=\Ass(\Oh_X/wJ)$ is a codimension~1 prime; the
lemma follows, since $\Ass(\Oh_X/J)=\Ass(w\Oh_X/wJ)\subset\Ass(\Oh_X/wJ)$.
Clearly,
 \[
 \Ass(\Oh_X/vI)\ \subset\ \Ass(\Oh_X/I)\cup\Ass(I/vI).
 \]
For any $P\in\Ass(I/vI)$, choose $x\in I$ with $P=(vI:x)=\Ann(\xbar\in I/vI)$.
One sees that
 \[
 \begin{cases}
 x\in \Oh_Xv \implies P\in\Ass(\Oh_X/I), \\
 x\notin \Oh_Xv \implies \hbox{$P\subset Q$ for some $Q\in\Ass(\Oh_X/v\Oh_X)$.}
 \end{cases}
 \]
Since every associated prime of $\Oh_X/v\Oh_X$ has codimension~1, this gives
 \[
 \Ass(\Oh_X/vI)\ \subset\ \Ass(\Oh_X/I)\cup\Ass(\Oh_X/v\Oh_X). \QED
 \]
\end{pf}

 \begin{theorem}[Kustin and Miller \cite{KM}] \label{thm!KM}
 The element $s\in \Oh_X[s]$ is a non-zerodivisor, and $\Oh_X[s]$ is a
Gorenstein ring.
 \end{theorem}

 \paragraph{Step 1} We first prove that
 \begin{equation}
 S\Oh_X[S]\cap(Sf_i-g_i)=S(Sf_i-g_i),
 \label{eq!coprime}
 \end{equation}
under the assumption that $s\colon I\to J$ is an isomorphism.

For suppose $b_i\in\Oh_X[S]$ are such that $\sum b_i(Sf_i-g_i)$ has no
constant term. Write $b_{i0}$ for the constant term in $b_i$, so that
$b_i-b_{i0}=Sb_i'$. Then $\sum b_{i0}g_i=0$. Since $s\colon f_i\mapsto
g_i$ is injective, also $\sum b_{i0}f_i=0$. Thus the constant terms in the
$b_i$ don't contribute to the sum $\sum b_i(Sf_i-g_i)$, which proves
(\ref{eq!coprime}).

The natural projection $\Oh_X[S]\onto\Oh_X$ takes $(Sf_i-g_i)\onto
J=(g_i)$, and (\ref{eq!coprime}) calculates the kernel. This gives the
following exact diagram:
 \[
 \renewcommand{\arraystretch}{1.5}
 \begin{matrix}
 0&\to&(Sf_i-g_i)&\xrightarrow{\ S\ }&(Sf_i-g_i)&\to&J&\to&0 \\
 && \bigcap && \bigcap && \bigcap \\
 0&\to&\Oh_X[S]&\xrightarrow{\ S\ }&\Oh_X[S]&\to&\Oh_X&\to&0 \\
 && \doubleheaddownarrow && \doubleheaddownarrow && \doubleheaddownarrow \\
 &&\Oh_X[s]&\xrightarrow{\ s\ }&\Oh_X[s]&\to&\Oh_X/J&\to&0
 \end{matrix}
 \]
The first part of the theorem follows by the Snake Lemma.

\paragraph {Step 2} To prove that $N$ is Cohen--Macaulay, recall that
 \[ 
 \depth M=\inf \big\{i \geq 0 \bigm| \Ext^i_{\Oh_X}(k,M)\ne0 \big\}
 \]
for $M$ a finite $\Oh_X$-module over a local ring $\Oh_X$ with residue
field $k=\Oh_X/\fm$ (see \cite{M}, Theorem~16.7). We have two exact
sequences 
 \begin{equation}
 \begin{gathered}
 0 \to I \to \Oh_X \to \Oh_X/I \to 0 \\[4pt]
 0 \to J \to \Oh_X \to \Oh_X/J \to 0. \label{eqn!ofJ}
 \end{gathered}
 \end{equation}
By assumption, $\Oh_X$ and $\Oh_X/I$ are Cohen--Macaulay, therefore 
\begin {gather*}
 \Ext^i_{\Oh_X}(k,\Oh_X)=0 \quad \hbox{ for } 0 \le i < n \\
\hbox{ and } \quad
 \Ext^i_{\Oh_X}(k,\Oh_X/I)=0 \quad \hbox{ for } 0 \le i < n-1,
\end{gather*}
where $n=\dim X$.
Thus
 \begin{equation}
 \Ext^i_{\Oh_X}(k,I)=0 \quad \hbox{ for } 0 \le i < n,
 \label{eq!dpth}
 \end{equation}
and the Ext long exact sequence of (\ref{eqn!ofJ}) gives also
 \[
 \Ext^i_{\Oh_X}(k,\Oh_X/J)=0 \quad \hbox{ for } 0 \le i < n-1.
 \] 
Therefore $\Oh_N=\Oh_X/J$ is Cohen--Macaulay.
 
\paragraph{Step 3} We prove that $\om_N\iso\Oh_N$ by running the argument
of Lemma~\ref{lem!inj} in reverse. Recall that $\sHom(I,\om_X)$ is generated
by two elements $i,s$, where $i$ is a given basis element of $\om_X$ viewed
as a submodule $\om_X\subset \sHom(I,\om_X)$, and $s$ is our isomorphism
$I\to J\subset\om_X$.

We write $j$ for the same basis element of $\om_X$ viewed as a submodule of
$\sHom(J,\om_X)$, and $t=s^{-1}\colon J\to I\subset\om_X$ for the inverse
isomorphism. Now $s\colon I\to J$ induces a dual isomorphism
 \[
 s^*\colon \sHom(J,\om_X) \to \sHom(I,\om_X),
 \]
which is defined by $s^*(\fie)(v)=\fie(s(v))$ for $\fie\colon J\to\om_X$.
By our definitions, clearly $s^*(j)=s$ and $s^*(t)=i$. Since $s^*$ is an
isomorphism, it follows that $\sHom(J,\om_X)$ is generated by $t$ and $j$.
Therefore the adjunction exact sequence
 \[
 0 \to \om_X \to \sHom(J,\om_X) \to \om_N \to 0
 \]
gives $\om_N=\Oh_N\tbar$. This completes the proof that $\Oh_N$ is
Gorenstein.

\paragraph{Alternative proof of Step~3} We worked out the above slick
proof of Step~3 by untangling the following essentially equivalent
argument, which may be more to the taste of some readers.

We set up the following exact commutative diagram:
 \[
 \renewcommand{\arraystretch}{1.5}
 \begin{matrix}
 0&\to&I&\xrightarrow{\ s\ }&\Oh_X&\to&\Oh_N&\to&0 \\
 && \bigcap && \bigcap \\
 0&\to&\Oh_X&\xrightarrow{\ s_2\ }&\sHom(I,\om_X)&\to&L&\to&0 \\
 && \doubleheaddownarrow && \doubleheaddownarrow \\
 &&\Oh_D&\xrightarrow{\ s_3\ }&\om_D \\
 \end{matrix}
 \]
The first column is just the definition of $\Oh_D$. The second column is
the identification of $\Oh_X$ with $\om_X$ composed with the adjunction
formula for $\om_D$.

The first row is the multiplication $s\colon I\to J$ composed with the
definition of $\Oh_N$. To make the first square commute, the map $s_2$
must be defined by
 \begin{equation}
 s_2(a)(b)=s(ab) \quad \hbox{for $a\in\Oh_X$ and $b\in I$.}
 \label{eq!s2}
 \end{equation}
We identify its cokernel $L$ below. The first two rows induce the map
$s_3$. Since $s_2$ takes $1\in\Oh_X$ to $s\in\sHom(I,\om_X)$, it follows
that $s_3$ takes $1\in\Oh_D$ to $\sbar\in\om_D$ as in Lemma~\ref{lem!inj},
and therefore $s_3$ is an isomorphism.

Now the second row is naturally identified with the adjunction sequence
 \[
 0\to\om_X\to\sHom(J,\om_X)\to\om_N\to0.
 \]
The point is just that $s\colon I\iso J$, and $s_2$ is the composite
 \[
 0\to\om_X\into\sHom(I,\om_X) \xrightarrow{s^*}\sHom(J,\om_X),
 \]
by its definition in (\ref{eq!s2}). The Snake Lemma now gives $\Oh_N\iso
L=\om_N$. Therefore, as before, $N$ is Gorenstein.

\medskip

In what follows, we prove that $\Oh_X[s]$ is Gorenstein: that is (see
\cite{M}, Definition on p.~145 and Theorem~18.2), its localisation
$(\Oh_X[s])_\fn$ at any maximal ideal $\fn$ of $\Oh_X[s]$ is Gorenstein.
Since $\fn\cap\Oh_X=\fp\in\Spec\Oh_X$, and localising $\Oh_X$ at $\fp$
preserves all the assumptions, we need only to consider $\fn$ lying over
$\fn\cap\Oh_X=\fm$.

\paragraph{Step 4 -- special case} If $s\in\fn$ we are done by Step~3:
$s$ is a regular element by the first part of the theorem and
$(\Oh_X[s])_\fn/(s)=\left(\Oh_N\right)_\fn$ is Gorenstein by Step~3. This
argument gives nothing if $s\notin\fn$: $s$ is a unit and the quotient by
$(s)$ is zero. (This was a small gap in the Nov 2000 preprint of this
paper.)

The same argument works if $s-a\in\fn$ for some $a\in\Oh_X$: we can
arrange that $s-a\colon I\to\Oh_X$ is injective (wiggle if necessary by
an element of $\fm$, as in Lemma~\ref{lem!inj}, (i)), then replace
$s\mapsto s-a$ in the construction by a coordinate change as in
Remark~\ref{rmk!expl}, (2). If $\Oh_X$ contains an algebraically closed
field $k$ that maps isomorphically to the residue field $\Oh_X/\fm$
(the main case in many applications), this completes the proof.

\paragraph{Step 5 -- general case} We use an extension of the residue
field $k=\Oh_X/\fm$ to reduce the general case to the case $s\in\fn$. We
need two facts.

 \begin{exc}\label{ex!O'} Let $(A,\fm,k)$ be a local ring, and $k\subset
L$ a finite extension of the residue field. Then there exists an extension
ring $A\subset B$ such that
 \begin{enumerate}
 \renewcommand{\labelenumi}{(\roman{enumi})}
 \item $B\iso A^{\oplus N}$ is a free $A$-module of rank $N=[L:k]$;
 \item $B$ is local with maximal ideal $\fm'=\fm\cdot B$ and $B/\fm'=L$.
 \end{enumerate}
[Hint: Do a primitive extension $k\subset k_1$ first, then induction on
$[L:k]$.]
 \end{exc}

 \begin{prop}\label{prop!W} Let $(A,\fm,k)$ be a local ring and $A\subset
B$ an extension ring that is a finite free $A$-module. Thus $B\iso
A^{\oplus N}$ and $B/\fm B\iso k^{\oplus N}$. In particular, $B$ is
semilocal, and its finitely many localisations $(B_i,\fn_i)$ also have
quotient rings $B_i/\fm B_i$ that are finite dimensional $k$-vector
spaces.
 \begin{enumerate}
 \renewcommand{\labelenumi}{(\roman{enumi})}
 \item $\depth A=\depth_{B_i}B_i$ for each $i$; in particular, $A$ is
Cohen--Macaulay $\iff$ $B$ is Cohen--Macaulay.
 \item $B$ Gorenstein $\implies$ $A$ Gorenstein. (In fact, it is enough
that one localisation $B_i$ is Gorenstein.)
 \item Let $A\subset B$ be as in \ref{ex!O'}, (i--ii). Then $A$ Gorenstein
$\implies$ $B$ Gorenstein.
 \end{enumerate}
 \end{prop}

 \begin{pf} This is a standard result in commutative algebra; see for
example Bruns and Herzog \cite{BH}, Theorem~1.2.16, p.~13, or Watanabe,
Ishikawa, Tachibana and Otsuka \cite{WITO}, or Matsu\-mura \cite{M},
Theorem~23.4. (We thank John Moody for explaining the argument of
\cite{WITO} to us.)

We sketch a direct proof to avoid these references. In (i), an
$A$-regular sequence remains $B_i$-regular, so $\depth A\le\depth B_i$ is
clear. The other way round, $\depth A=0$ means that $\fm\in\Ass A$, that
is, $A$ contains a copy of $A/\fm$. Then, by flatness, each $B_i$
contains a copy of $B_i/\fm B_i$. This is an Artinian local ring, so
$\{\fn_i\}=\Ass B_i/\fm B_i\subset\Ass B_i$, and therefore
$\depth_{B_i}B_i=0$. Thus $\depth B_i>0$ implies also $\depth A>0$. (i)
follows by the usual induction on $\depth A$.

After dividing out by a maximal $A$-regular sequence
$x_1,\dots,x_n\in\fm$, in (ii) and (iii) we can assume that $A$ and
$B_i$ are Artinian, and everything comes down to estimates on the
dimension of socles. For (ii), suppose that $(A/\fm)^{\oplus
a}\subset\socle A$. As above, $(B_i/\fm B_i)^{\oplus a}\subset B_i$ by
flatness, and $B_i/\fm B_i$ contains at least one copy of $B_i/\fn_i$.
Thus $(B_i/\fn_i)^{\oplus a}\subset\socle B_i$. Thus $B_i$ Gorenstein
implies $a=1$ and $A$ is Gorenstein. For (iii), $\dim\Hom_A(k,A)=1$ and
$B\iso A^{\oplus N}$ gives $\dim\Hom_A(k,B)=N$ with $N=[L:k]$. On the
other hand, if $L^{\oplus a}\subset\socle B$ then $\dim\Hom_A(k,B)\ge
a[L:k]$; therefore $a\le1$, and $B$ is Gorenstein. \QED \end{pf}

We return to the proof of Theorem \ref{thm!KM}. The maximal ideal of the
polynomial ring $\Oh_X[S]$ lying over $\fn\subset\Oh_X[s]$ is of the form
$(\fm,F)$ for some monic polynomial $F\in\Oh_X[S]$ whose reduction
$\Fbar\in k[S]$ remains irreducible over $k$. Write $k\subset L$ for a
splitting field of $\Fbar$, so that $\Fbar=\prod_i(S-\al_i)$ with
$\al_i\in L$, and repeated factors if $\Fbar$ is inseparable. Let
$\Oh_X\subset\Oh'$ be a ring extension as in \ref{ex!O'}, (i--ii). We
write $I'=I\Oh'$ and $J'=J\Oh'$, and extend $s$ to an isomorphism $I'\iso
J'$ given by the same formula $f_i\mapsto g_i$. The extended rings $\Oh'$
and $\Oh'/I'$ are Gorenstein by Proposition~\ref{prop!W}, (iii). Thus
 \[
 \Oh'[s]=\Oh'[S]/(f_iS-g_i)=\Oh_X[s]\tensor\Oh'
 \]
is the unprojection ring constructed from the local ring $\Oh'$ and the
ideal $I'$.

Now every maximal ideal of $\Oh'[s]$ that contains $\fn$ is of the
form $\fn_i'=(\fm,s-a_i)\Oh'$, and has residue field $L$, where
$a_i\in\Oh'$ reduces to $\al_i\in L$. Thus each localisation
$(\Oh'[s])_{\fn_i'}$ is Gorenstein by the argument of Step~4, and
Proposition~\ref{prop!W}, (ii) gives $(\Oh_X[s])_\fn$ Gorenstein.

This completes the proof of Theorem~\ref{thm!KM}.  \QED

\subsection*{Kustin and Miller's argument}
We paraphrase the argument of \cite{KM} for completeness. In addition to
our usual assumptions, suppose that everything is contained in an ambient
local scheme $A$, with $\Oh_X$ and $\Oh_D$ of finite projective dimension
over $\Oh_A$ (for example, because $A$ is regular). Suppose that
$\codim_A X=d$ and $\codim_A D=d+1$. We write out free resolutions
$L\lodot\to\Oh_X$ and $M\lodot\to\Oh_D$ over $\Oh_A$. Then the usual
properties of resolutions give a map of complexes
 \begin{equation}
 \renewcommand{\arraystretch}{1.5}
 \begin{array}{cccccccccc}
 & 0\to& L_d & \cdots & L_1 &\to& \Oh_A &\to&\Oh_X \\
 &&\big\downarrow && \big\downarrow &&\big\Vert && \doubleheaddownarrow \\
 0\to M_{d+1} & \to & M_d & \cdots & M_1 &\to& \Oh_A &\to& \Oh_D
 \end{array}
 \label{eq!cx_map}
 \end{equation}
Suppose that the ideal of $D$ in $\Oh_A$ is generated by $k$ elements
$f_1,\dots, f_k$, so that $M_1=k\Oh_A$.

We identify $L_d=M_{d+1}=\Oh_A$ and $M_d=M_1^\vee=k\Oh_A$ by Gorenstein
symmetry. Then the tail-end of the complexes gives
 \begin{equation}
 \renewcommand{\arraystretch}{1.5}
 \begin{array}{ccc}
  & L_d=\Oh_A&\to\cdots \\
  & \big\downarrow & \kern-8mm(g_1,\dots,g_k) \\
 0 \to M_{d+1}=\Oh_A \xrightarrow{f_1,\dots,f_k} & M_d=k\Oh_A &\to\cdots
 \end{array}
 \label{eq!tail}
 \end{equation}

As we have done, Kustin and Miller introduce a new indeterminate $S$, and
write out new equations $Sf_i=g_i$. This gives a new ambient space
$\aff^1_A=\Spec\Oh_A[S]$ and a new ideal $I_Y=(I_X,Sf_i-g_i)$. They introduce
a new complex by glueing together $L\lodot\tensor \Oh_A[S]$ and
$M\lodot\tensor \Oh_A[S]$, and prove it is a resolution of $I_Y$ by
arguments based on the Buchsbaum--Eisenbud criterion \cite{BE}.

To check that their construction is the same as ours, take the dual of
(\ref{eq!tail}), note that $\om_X=\coker \bigl\{L_{d-1}^\vee\to
L_d^\vee\bigr\}$, and identify the ideal of $D$ in $A$ with
$\coker\big\{(f_i)\colon M_d^\vee\to M_{d+1}^\vee\big\}$. One shows that
the dual diagram induces a map $f_i\mapsto g_i$ from
$(f_i)=I_{X,D}\subset\Oh_X$ to $\om_X$ that provides the second
generator of $\sHom(I,\om_X)$ as in Lemma~\ref{lem!inj}. See \cite{P1},
Section~3 for details.

The advantage of their method is that it gives in theory the complex
resolving the new ideal. On the other hand, while it is trivial to say
that the map of complexes (\ref{eq!cx_map}) exists, it is hard to
calculate, except in the simplest examples (compare \ref{ssec!lit}
below); some cases are worked out in Papadakis \cite{P}--\cite{P1} and
\cite{Ki}. Our construction identifies the final and most important
vertical map in (\ref{eq!cx_map}) as a Poincar\'e residue, and thus as a
``determinant'', which we can hope to calculate birationally without
knowing the finer details of the ``matrix'' that gave rise to it. (From
a philosophical point of view, this is the whole point of the canonical
class!)

\section{Applications}\label{sec:appl}
\subsection{The affine case} Consider first the geometry of the affine
graph: the morphism $\pi\colon Y\to X$ is the graph of $s$. The locus
$D\setminus N$, where $s$ has a pole, disappears off ``to infinity'' on
$Y$, whereas $N\setminus D$ becomes the principal divisor $s=0$. The
intersection $D\cap N$ is the locus of indeterminacy of $s$, and $Y$
contains an affine line bundle over it. 

\subsection{Example: nodal curve}\label{ssec!nodal}
 The example $X=$ nodal curve, $D=$ reduced node is very instructive;
set $X:(x^2-y^2=0)$ and $D:(x=y=0)$. Then $s=x/y$ is an automorphism of
$I=J=\fm$, and $Y\to X$ is an affine blowup, with an exceptional
$\aff^1$ over the node. (The affine line in question is the
projectivised of the 2-dimensional vector space
$\Hom_{\Oh_X}(\fm,\fm)\tensor k$, with the identity element deleted.)
This example shows that Remark~\ref{rmk!expl}, (4) does not hold
without the assumption that $D$ and $N$ are Weil divisors.

\subsection{Simplest example}\label{ssec!sim}
We discuss a case that has many consequences in birational geometry, even
though the algebra itself is very simple. Consider the generic equations
 \begin{equation}
 X:(Bx-Ay=0)\quad\hbox{and}\quad D:(x=y=0)
 \label{eq!Ax-By}
 \end{equation}
defining a hypersurface $X$ containing a codimension~2 complete intersection
$D$ in some as yet unspecified ambient space. The unprojection variable is
 \begin{equation}
 s=\frac{A}{x}=\frac{B}{y}\,.
 \label{eq!deg1}
 \end{equation}
We can view $s$ as a rational function on $X$, or as an isomorphism from
$(x,y)$ to $(A,B)$ in $\Oh_X$. The unprojection is the codimension two
complete intersection $Y:(sx=A,sy=B)$.

For example, take $\PP^3$ as ambient space, with $x,y$ linear forms defining
a line $D$, and $A,B$ general quadratic forms. See \ref{ssec!prj} for the
standard trick to make our local construction work also in the projective
set-up. Then $s$ has degree~1 from (\ref{eq!deg1}), and the equations
describe the contraction of a line on a nonsingular cubic surface to the
point $P_s=(0:0:0:0:1)\in\PP^4$ on a del Pezzo surface of degree~4. It is the
inverse of the linear projection $Y\broken X$ from $P_s$, eliminating $s$.
But the equations are of course much more general. The only assumptions are
that $x,y$ and $Bx-Ay$ are regular sequences in the ambient space. For
example, if $A,B$ vanish along $D$, so that $X$ is singular there, then
$Y$ contains the plane $x=y=0$ as an exceptional component lying over $D$
(as in \ref{ssec!nodal}). Note that, in any case, $Y$ has codimension~2 and
is nonsingular at $P$.

The same rather trivial algebra lies behind the quadratic involutions of
Fano \hbox{3-folds} constructed in Corti, Pukhlikov and Reid \cite{CPR},
4.4--4.9. For example, consider the general weighted hypersurface of
degree~5
 \[
 X_5:(x_0y^2+a_3y+b_5=0) \subset \PP(1,1,1,1,2),
 \]
with coordinates $x_0,\dots,x_3,y$. The coordinate point $P_y=(0:\cdots:1)$
is a Veronese cone singularity $\frac12(1,1,1)$. The anticanonical model of
the blowup of $P_y$ is obtained by eliminating $y$ and adjoining $z=x_0y$
instead, thus passing to the hypersurface
 \[
 Z_6: (z^2+a_3z+x_0b_5=0) \subset \PP(1,1,1,1,3).
 \]
The 3-fold $Z_6$ contains the plane $x_0=z=0$, the exceptional $\PP^2$ of
the blowup. Writing its equation as $z(z+a_3)+x_0b_5$ gives
$y=\frac{z}{x_0}=-\frac{b_5}{z+a_3}$, and puts the birational relation
between $X_5$ and $Z_6$ into the generic form
(\ref{eq!Ax-By}--\ref{eq!deg1}). In fact $Z_6$ is the ``midpoint'' of the
construction of the birational involution of $X_5$. The construction
continues by setting $y'=\frac{z+a_3}{x_0}=-\frac{b_5}{z}$, thus
unprojecting a different plane $x_0=z+a_3=0$. For details, consult
\cite{CPR}, 4.4--4.9. See Corti and Mella \cite{CM} for a related use of
the same algebra, to somewhat surprising effect; these and many further
examples are treated at more length in Papadakis \cite{P}--\cite{P1} and
Reid \cite{Ki}.

\subsection{Projective case}\label{ssec!prj}
The cases of \ref{ssec!sim} are typical of our applications of unprojection
to biregular and birational geometry. Although we developed the theory for
local rings in Section~1, it applies at once to projective varieties
(schemes) via the standard philo\-sophy of Zariski and Serre summarised in
the slogan ``graded is a particular case of local'' (the coherent half of
Grothendieck's ``Lefschetz principle'', see Grothendieck \cite{G}, esp.\
Chapters~I--V). We sketch briefly what we need.

Our graded rings $R$ are graded in positive degrees, with $R_0=k$ a
field, and $R$ finitely generated over $k$. The associated local ring
$R_\fm$ is $R$ localised at the maximal ideal $\fm=\bigoplus_{n>0}R_n$.
The principle says that coherent co\-homo\-logy of sheaves on $X$ can be
treated in terms of local co\-homo\-logy $H_\fm^i(M)$ of modules over
$R_\fm$ at $\fm$. In particular, Cohen--Macaulay and Gorenstein have
equivalent treatments in terms of the geometry of the projective scheme
$X$ or the local cohomology of $R_\fm$. Geometrically, $\Spec R$ is the
affine cone over $X=\Proj R$, and we localise at the origin $(0)=V(\fm)$;
this replaces the projective variety by a small neighbourhood of the
vertex of its affine cone (together with its grading).

Let $I\subset R$ be a graded ideal of codimension~1. We suppose that $R$
and $R/I$ are Gorenstein and write $X=\Proj R$ and $D=\Proj R/I$. Then
$D\subset X$ are projectively Gorenstein schemes. Write $\om_X=\Oh_X(k_X)$
and $\om_D=\Oh_D(k_D)$, and assume that $k_X>k_D$ (see
Remark~\ref{rem!kX>kD}). The construction of Section~1 gives a rational
section $s$ of $\Oh_X(k_X-k_D)$ that defines an isomorphism $I\to J$ of
ideals of $R$, but with a Serre twist by $k_X-k_D$. (For example, in
\ref{ssec!sim}, $s\colon(x,y)\mapsto(A,B)$ has degree~1.) It is naturally
an isomorphism of sheaves $s\colon I\iso J(k_X-k_D)$. As before, write
$R[s]:=R[S]/(Sf_i-g_i)$ for the unprojection ring and $Y=\Proj R[s]$.

If $R=k[x_1,\dots,x_k]/(I_X)$ is generated by elements $x_i$ with $\deg
x_i=a_i$ then $R[s]=k[x_1,\dots,x_k,s]/(I_Y)$ is generated by $x_i$ together
with $s$, of degree $\deg s=k_X-k_D$. The projective scheme $Y$ contains the
distinguished point $P_s=(0:\dots:0:1)$. If $X$ is variety (that is,
reduced and irreducible), and $D$ a Weil divisor, then $D\cap N$ does not
contain any prime divisors, so that the inclusion $R\subset R[s]$ defines a
birational map or {\em unprojection} $X\broken Y$ contracting $D$ to the
point $P_s$. This is the striking difference from the affine case, where
$D\setminus N$ disappeared ``off to infinity''. The inverse rational map
$Y\broken X$ is the projection from $P_s$, and corresponds to eliminating
$s$. It blows up $P_s$ to a divisor $E\subset\wY$, then defines a morphism
$\wY\to Y$ taking $E$ birationally to $D$.

 \begin{rmk}\label{rem!kX>kD}
 We need the assumption $\deg s=k_X-k_D>0$ in order that the unprojection
ring is still graded in positive degrees. It is interesting to note that it
also has the effect of making $D$ have ``negative self-intersection''. If
$D$ is a Cartier divisor on $X$ then of course $\Oh_D(-D)=\Oh_D(k_X-k_D)$.
However, even if $D\subset\Sing X$, so that it is not even a Weil divisor,
the difference $(\om_X)\rest{D}\tensor\om_D\1$ provides a usable notion of
self-intersection class of $D$. For more discussion of this idea, compare
Catanese, Franciosi, Hulek and Reid \cite{CFHR}, 3.1--3.2 and \cite{Ki},
Section~5.3.
 \end{rmk}

 \begin{rmk} The case $\dim X=1$ leads to an apparent paradox, kindly
brought to our attention by Nikos Tziolas: Theorem~\ref{thm!KM} applied to
a point on a canonical curve $P\in C\subset\PP^{g-1}$ seems at first sight
to give a canonical form with a single pole at $P$. The source of the
misunderstanding is that if $\dim D=0$ then $D$ (up to isomorphism) does
not determine $k_D$; the solution is to keep track of the graded ring over
$D$, which does. In fact $P=\Proj k[x]$, with $\deg x=1$, so that
$k_P=-1$, and $k_C-k_P=2$. Our construction in this case unprojects $C$ to
a curve $C\subset\PP(1^g,2)$. This curve is projectively Gorenstein, with
$\Oh_C(1)$ corresponding to the $\Q$-Cartier divisor $K_C+\frac12P$. It is
a nonsingular curve passing through the $\Z/2$ quotient singular
$P_s=(0,\dots,0,1)$, and should be viewed as an orbifold at that point,
with orbifold canonical class $K_C+\frac12P$.
 \end{rmk}

 \begin{ex} The humane treatment of the combinatorics of the graded
ring $R(X,\Oh_X(1))$ is in terms of its Poincar\'e series
$P_X(t)=\sum_{n\ge0} P_n(X)t^n$, where $P_n(X)=h^0(X,\Oh_X(n))$ (see
for example Alt{\i}nok \cite{A1}). If $X,D$ and $Y$ are as above, prove
that
 \[
 P_Y(t) = P_X(t) + \frac{t^k}{1-t^k}\, P_D(t),\quad
 \hbox{where } k=\deg s = k_X-k_D.
 \]
 \end{ex}

 \subsection{Cases already in the literature}\label{ssec!lit}
 \cite{KM} and \cite{KM1} contain many examples of unprojections in
commutative algebra. The following case (already in \cite{KM0}) is probably
the simplest: let $A=a_{ij}$ be a generic $3\times4$ matrix and
$x=(x_1,\dots,x_4)$ a $4\times1$ column vector. Write $q_i=\sum a_{ij} x_j$;
these are 3 generic linear combinations of a regular sequence of length~4.
Then $D:(x_j=0)$ is a codimension~4 complete intersection, and $X:(q_i=0)$
a codimension~3 complete intersection containing $D$. The unprojection of
$D$ in $X$ consists of introducing a new indeterminate $s$ and writing out
the equations
 \begin{equation}
 Ax = 0, \quad sx = \bigwedge\nolimits^3 A.
 \label{eq!7x12}
 \end{equation}
That is, $sx_i$ equals $\pm$ the $3\times3$ minor $\det A_i$ obtained by
deleting the $i$th row of $A$. These equations are of course familiar from
Cramer's rule. In this case, the complexes $L\lodot$ and $M\lodot$ of
(\ref{eq!cx_map}) are Koszul complexes, and the vertical arrows are
successive wedges of $A$.

Equations (\ref{eq!7x12}) give one of the simplest formats of Gorenstein
rings in codimension~4, with a $7\times12$ resolution (7 relations, 12
syzygies). For example: $X\subset\PP^6$ a K3 surface given as a complete
intersection of 3 quadrics, containing a line $D$. The unprojection
contracts $D$ to an ordinary node of a K3 surface $Y\subset \PP(1^7,2)$
with ideal defined by (\ref{eq!7x12}).

It is quite curious that the format (\ref{eq!7x12}) occurs very rarely in
nature, possibly because it does not have a grading with all generators of
degree~1. For example, the canonical ring $R(X,K_X)$ of a regular surface of
general type with $p_g=5,K^2=11$ has Hilbert series $P_X(t)$ given by
 \begin{align*}
 P_X(t) & \ = \ \frac{1+(p_g-3)t+(K^2-2p_g+4)t^2+(p_g-3)t^3+t^4}{(1-t)^3} \\
 & \ = \ \frac{1-6t^3-t^4+12t^5-t^6-6t^7-t^{10}}{(1-t)^5(1-t^2)^2},
 \end{align*}
so a first hope is that $X\subset\PP(1^5,2^2)$ could be defined by 6
relations in degree~3, and 1 relation in degree~4, yoked by 12 syzygies in
degree~5. In fact, for geometric reasons, the two generators $y_1,y_2$ of
degree~2 must occur in 3 degree~4 relations $y_1^2=\cdots$, etc., so that
the resolution must be at least $9\times16$.

In contrast, rings with $9\times16$ resolutions are ubiquitous, and we
know a couple of hundred cases. Most of the families of K3s and Fano
3-folds in codimension~4 studied in Alt{\i}nok \cite{A} (more than a
hundred of them; see also \cite{R2}) are unprojections in weighted
projective space. Their rings almost invariably have $9\times16$
resolutions. There is a vaguely formulated (but strongly documented)
conjecture that rings with $9\times16$ resolutions fall into 2 standard
families called Tom and Jerry, bearing a family resemblance to the
Segre embeddings of $\PP^2\times\PP^2$ and
$\PP^1\times\PP^1\times\PP^1$. See \cite{P} and \cite{P1} for details
of the families as unprojections and compare \cite{Ki}, Section~8.

\subsection{Gorenstein projections and discrepancy} When we project a
del Pezzo surface $Y$ of degree $d\ge2$ to another $X$ of degree $d-1$,
the original surface satisfies $K_Y=-H$, the blowup $\wY\to Y$
introduces a $-1$-curve $E$ of discrepancy 1, so that $K_{\wY}=-\wH+E$,
and the projection morphism $\wY\to X$ is given by $\wH-E=-K_{\wY}$.
Thus the condition that both $X$ and $Y$ are projectively Gorenstein
(here, anti\-canonically polarised) involves an implicit discrepancy
calculation. We avoid a formal definition of Gorenstein projection
(except in the simplest case given below), because the notion makes
sense at different levels of generality. We only explain briefly the
point about discrepancy for the reader versed in 3-folds.

Consider a projectively Gorenstein variety $Y,H$ with $K_Y=kH$, a point
blowup $\si\colon \wY\to Y$ extracting a divisor $E$ with discrepancy $a$,
with a new divisor $\wH=H-mE$ that defines a small contraction $\fie\colon
\wY\to X$ with $\fie(E)=D\subset X$. Suppose that also $K_X=kH_X$; then
because $\wY\to X$ is small, also $K_{\wY}=k\wH$, and
 \begin{equation}
 K_{\wY} = K_Y+aE, \quad \wH = H-mE \implies a+km = 0.
 \label{eq!discr}
 \end{equation}
For Fano 3-folds, $k=-1$, and the basic operation in \cite{CPR} is the
Kawamata blowup of a quotient singularity $\frac{1}{r}(1,a,r-a)$, which
satisfies $a=m=\frac{1}{r}$. For $Y$ a K3 surface with Du Val
singularities, $k=0$, and any crepant blowup (that is, with $a=0$) may
lead to a Gorenstein projection (and does so provided that $H-mE$ is still
nef and big). Another interesting case is a canonically embedded regular
surface $Y$ with (say) a simple elliptic singularity. Then $k=1$, $a=-1$,
and the linear projection may be a Gorenstein projection; this type of
blowup decreases $p_aY$ by 1, and $K_Y^2$ by the degree of the
singularity, and appears in various forms in many constructions of
algebraic surfaces of general type.

We give a temporary definition to serve as a converse to
Theorem~\ref{thm!KM}. Say that $Y\broken X$ is a {\em Gorenstein
projection of Type~I\/} if $\si\colon\wY\to Y$ is a point extraction
with exceptional divisor $E$, and $\fie\colon\wY\to X$ a small
birational morphism taking $E$ birationally to $D\subset X$, such that
$X$ and $D$ are both projectively Gorenstein. In this case $X,D$
satisfies the assumptions of Theorem~\ref{thm!KM}, and the unprojection
leads back to $Y$. Necessary conditions for this to hold are the
restriction (\ref{eq!discr}) on the discrepancy of the blowup, and the
surjectivity of the restriction maps $H^0(\wY,n\wH)\to H^0(E,n\wH\rest
E)$ (needed so that $D$ is projectively normal).

\subsection{More general $D$} Finally, our construction has natural
generalisations: the exceptional divisor $D\subset X$ of a Gorenstein
projection is not restricted to being projectively Gorenstein. For
example, Fano's construction of his 3-folds $V_{2g-2}\subset\PP^{g+1}$ for
$g\ge6$ involves projections $Y\broken X$ from a line, giving rise to the
exceptional locus $D\subset X$ a cubic scroll. Alessio Corti observed that
the inverse map can be treated by a generalisation of our method: to
unproject $D$, we need to take two generators of $\sHom(I_D,\om_X)$ of
degree~1, that map down to a basis of $\om_D(2)$, the ruling of the
scroll. (He also carried out in some detail the algebraic calculation of
the anticanonical ring of $Y$ in these terms.)

In the weighted projective case, it can happen that
$\om_D\iso\Oh_D(k_D)$,  but $D$ is not projectively normal. For
example, consider the following method for constructing Fano 3-folds:
let $D$ be a weighted projective plane $\PP(1,b_1,b_2)$ (for some small
values of $b_1,b_2$), embedded as a projectively nonnormal surface in
some $\PP^4(1,a_1,a_2,a_3,a_4)$. A Fano hypersurface $X_d$ with mild
singularities can contain $D$, and then $D$ can be contracted in $X$ to a
terminal quotient singularity $\frac{1}{b_1+b_2}(1,b_1,b_2)$ in a Fano
3-fold $Y$ by an unprojection.

Thus $D$ can sometimes be unprojected even though it is not Gorenstein, for
example because it is not normal. Arguing locally as in Section~1, suppose
that $D$ is a divisor in a normal Gorenstein variety $X$, that $D$ is
nonnormal, but its normalisation map $\wD\to D$ is an isomorphism in
codimension~1, and such that $\wD$ is Gorenstein. For example, $D$ could be
$\aff^2$ with the subscheme $(x^2=y=0)$ pinched to a point, or $\aff^3$
with the $x$-axis glued to itself by $x\mapsto-x$. Since
$\om_D=\om_{\wD}$, it is generated by a single element $\sbar$ as a module
over the normalisation $\Oh_{\wD}$. Thus $\sHom(I_D,\om_X)$ contains a
generator $s$ mapping to the basis $\sbar\in\om_{\wD}$. From the
birational point of view, $s$ has a pole along $D$, so its graph sends $D$
off to a point at infinity. Thus adjoining $s$ unprojects $D$ in a variety
$Y_0$. We believe that the normalisation $Y\to Y_0$ is a Gorenstein
variety (see \cite{Ki}, Problem~9.2). In this case, the algebraic
presentation of the unprojection ring $\Oh_Y$ is of course much more
complicated, and only worked out in special cases. The normalisation
amounts to taking a few extra generators $s_0=s,s_1,\dots,s_k\in
\sHom(I_D,\om_X)$, corresponding to generators of $\om_D$ over $\Oh_D$
(rather than $\Oh_{\wD}$), and adjoining these to $\Oh_X$. The case we
understand is discussed in \cite{CPR}, 7.3 and \cite{Ki}, Section~9; see
also the corresponding discussion in \cite{CM} and \cite{T4}. Although
the unprojection ring $\Oh_X[s_0,\dots,s_k]$ in this case certainly has
relations of the form $s_is_j=\cdots$ that are quadratic in the $s_i$,
it is still uniquely determined by its linear relations, essentially
coming from the presentation of $\om_D$ over $\Oh_D$.

\bigskip
\noindent
Stavros Papadakis, \\
Math Inst., Univ. of Warwick,\\
Coventry CV4 7AL, England\\
e-mail: spapad@maths.warwick.ac.uk

\medskip
\noindent
and Karaiskaki 11, Agia Paraskevi, \\
GR 153 41, Attiki, Greece \\
e-mail: stavrospapadakis@hotmail.com

\medskip
\noindent
Miles Reid,\\
Math Inst., Univ. of Warwick,\\
Coventry CV4 7AL, England\\
e-mail: miles@maths.warwick.ac.uk \\
web: www.maths.warwick.ac.uk/$\!\sim$miles

\end{document}